\documentclass{amsart}[12pt]
 \usepackage{geometry}                % See geometry.pdf to learn the layout options. There are lots.
 \usepackage{graphicx}
 \usepackage{amssymb}
 \usepackage{epstopdf}
\usepackage{amsmath}
\usepackage{amsfonts}
\usepackage{mathtools}
\usepackage{amsthm}
\usepackage{indentfirst}

\newtheorem{theorem}{Theorem}[section]
\newtheorem{lemma}[theorem]{Lemma}

\newtheorem{definition}[theorem]{Definition}

\newtheorem{proposition}[theorem]{Proposition}

\newtheorem{remark}[theorem]{Remark}

\newtheorem{corollary}[theorem]{Corollary}

\begin{document}

%end of theorem environments

%\begin{center}

\title[Nil-Good and Nil-Good Clean Matrix Rings]{Nil-Good and Nil-Good Clean Matrix Rings}
% \title[Sums of units and nilpotents in rings]{Sums of units and nilpotents in rings}
\author[Block Gorman]{Alexi Block Gorman}
\author[Shiao]{Wing Yan Shiao}

\address{Department of Mathematics, Wellesley College, Wellesley, MA 02481}
\email{ablockgo@wellesley.edu}
\address{Department of Mathematics, Wellesley College, Wellesley, MA 02481}
\email{wshiao@wellesley.edu}

%\end{center}

\maketitle

\begin{abstract}
The notion of clean rings and 2-good rings have many variations, and have been widely studied.  We provide a few results about two new variations of these concepts and discuss the theory that ties these variations to objects and properties of interest to noncommutative algebraists. A ring is called ``nil-good'' if each element in the ring is the sum of a nilpotent element and either a unit or zero.  We establish that the ring of endomorphisms of a module over a division is nil-good, as well as some basic consequences.  We then define a new property we call ``nil-good clean,'' the condition that an element of a ring is the sum of a nilpotent, an idempotent, and a unit.  We explore the interplay between these properties and the notion of clean rings.
\end{abstract}

\section{Introduction}

In 1977, W. K. Nicholson defined a ring $R$ to be clean if for every $a \in R$ there is $u$ a unit in $R$ and $e$ an idempotent in $R$ such that $a=u+e$ \cite{Nic77}.  The interest in the clean property of rings stems from its close connection to exchange rings, since clean is a concise property that implies exchange.  Properties of rings related to the clean and exchange properties have been largely expanded and researched, and some generalizations closely relate to other properties of interest to algebraists.

The study of rings generated by their units dates back to the 1950's, when it was proved that any endomorphism of a module over a division ring is equal to the sum of two units, unless the dimension of the module is 1 or the division ring is $\mathbb{F}_2$, as established in \cite{Wol53} and \cite{Zel54}.  This motivated algebraists to make extensive study of rings generated by their units.  Later, Peter V{\'a}mos defined an element $a$ in $R$ to be $2$-$good$ if it can be expressed as the sum of two units in $R$, and defined a ring $R$ to be $2$-$good$ if every element in $R$ is $2$-$good$ \cite{Vam05}.  In general, a ring is $n$-$good$ if each element can be written as the sum of $n$ units, and these properties have distinct applications from those of clean, but have also led to a diverse line of inquiry in ring theory.  P. Danchev defined a property in \cite{Dan15} related to $2$-$good$ in the following way: an element $a$ in $R$ is nil-good if $a=n+u$ where $n$ is a nilpotent element of $R$ and $u$ is either $0$ or a unit in $R$.  The ring $R$ is called nil-good if every element of $R$ is nil-good.

In this paper, we prove that if $R$ is a division ring, then $\mathbb{M}_n(R)$ is nil-good for all $n \in \mathbb{N}$.  We then establish some basic properties of nil-good rings in general.  We extend these results to specifically characterize local rings and artinian rings that are nil-good.

We then relate this property to clean rings in a new way by defining the property nil-good clean.  We say a ring $R$ is nil-good clean if for all $r \in R$ there is a unit $u$, a nilpotent $n$ and an idempotent $e$ in $R$ such that $r=u+n+e$.  This property holds for all clean and all nil-good rings, yet we show it includes a larger class of rings than only those that satisfy one property or the other.  Understanding how the nil-good clean property generalizes to include rings that are neither nil-good nor clean may reveal more about the interaction of those two properties within unital rings.   We provide an example of a nil-good clean ring that is not exchange, and therefore not clean.  The example has properties similar to the ring provided by G. M. Bergman in \cite{Han77} of a nonclean exchange ring.

Throughout this paper rings are associative with unity.  We denote the Jacobson radical $J(R)$ for a ring $R$ and write $\mathbb{M}_n(R)$ for the ring of $n \times n$  matrices over $R$.

\section{Matrices over a Field}

We first prove that $\mathbb{M}_n(R)$ is a nil-good ring when $R$ is a field for illustrative purposes. The linear algebra over $\mathbb{M}_n(R)$ required in the case where R is a field is more accessible, and the subsequent proof in the case where $R$ is a division ring is more concise and intuitive as a result.  We can write a nil-good decomposition for any element of $\mathbb{M}_n(k)$ where $k$ is a field by putting all matrices in rational canonical form.

\begin{theorem}{For all $n \geq 1$ the ring $\mathbb{M}_n(k)$ is nil-good if $k$ is a field.}
\begin{proof}

Note that suitable rearrangement of the basis elements allows us to rearrange the companion matrix blocks in a matrix's rational canonical form without altering the nilpotence or invertibility of that matrix.

If the coefficients of any matrix are all zero then its rational canonical form is exclusively zero except possibly on the subdiagonal, making it nilpotent, in which case we let $U$ be the zero matrix, and let $N$ be the matrix in question.
Similarly, if the matrix $A$ we wish to decompose is the $n \times n$ zero matrix we let both $U$ and $N$ be the zero matrix.  If the matrix $A$ in question is a unit, we let $N$ be the zero matrix, and let $U$ be the original matrix, the rational canonical form of $A$.

Now suppose $A$ is a non-nilpotent, non-unit matrix.  Then, it will have zero as its $-c_0$ coefficient for some companion matrix block, as well as some non-zero $-c_0$ coefficients.  We may choose an arrangement of the companion matrix basis that allows us to place the companion matrix blocks corresponding to the nonzero first coefficient in the upper left corner, and to place the blocks for which the $-c_0$ coefficient is zero in the lower right hand corner, ordered amongst themselves by size of block. \\
Thus, we consider a matrix of the form\\
 %\bigskip
 $$\begin{bmatrix}
            C_{g_1}  & 0  &\cdots & 0       \\
            0 & C_{g_2}    & \ddots  & \vdots  \\
	    \vdots   &\ddots    & \ddots   &  0  \\
            0           & \cdots         &0    &C_{g_k}
          \end{bmatrix}$$ \\
where each $C_{g_i}$ is a companion matrix that has some nonzero element for all $i<j$ for some $j>1$.  We call an $r \times r$ block consisting of zeros everywhere except the subdiagonal, which consists entirely of ones, an $N_r$ block. Therefore, each $C_{g_i}$ for $i<k$ has some nonzero coefficient $-c_i$ or is an $N_r$ block of size $2 \times 2$ or greater.

If the size $m \times m$ companion matrix block $C_{g_i}$ of $A$ is invertible, then we let the corresponding
diagonal block of $U$ be $C_{g_i}$ and the corresponding diagonal block of $N$ be the $m \times m$ zero block.

If $A$ contains any companion matrix block with zero as the $-c_0$ coefficient, then in the unit $U$ of its decomposition we add a $-1$ in the entry corresponding to the $-c_0$-coefficient of the companion matrix block, so that the block becomes invertible.   Correspondingly, we place a $1$ in the entry of the nilpotent matrix corresponding to the entry of the $-c_0$ coefficient of that companion matrix block. Otherwise we leave the corresponding block in the nilpotent summand entirely zero, so that the direct sum of these blocks will be the direct sum of nilpotents, making the overall matrix nilpotent as well.

Note that the method described works for $N_r$ blocks as well as those with nonzero coefficents.  In general, the decomposition of the companion matrix block that has zero as its $-c_0$ coefficient looks like:
$$\begin{bmatrix}
            0        & \cdots   & \cdots   & 0  \\
	     1  &\ddots  & \vdots   &  -c_1  \\
             0   &\ddots & 0    &\vdots    \\
            \vdots   &\ddots      &1   &  -c_{m-1}
          \end{bmatrix} = \begin{bmatrix}
 		0        & \cdots   &  \cdots & -1  \\
	     1  &\ddots  & \vdots   &  -c_1  \\
             0   &\ddots & 0    &\vdots    \\
            \vdots   &\ddots      &1   &  -c_{m-1}
          \end{bmatrix} +  \begin{bmatrix}
		0        & \cdots   & 0   & 1  \\
	     \vdots &\ddots  &    &  0  \\
             \vdots   &\ddots & \ddots    &\vdots    \\
            0   & \cdots   & \cdots  &  0
          \end{bmatrix}$$

Since these matrices will be useful later, the first matrix on the right side of the equation we will denote $C^*_{g_i}$ and the second matrix on the right we will denote $N^*$.  Observe that $C^*_{g_i}$  is invertible and $N^*$ is nilpotent.

Whether the coefficients $-c_1$ through $-c_{m-1}$ are zero or nonzero does not affect the validity of the decomposition, since the addition of a one in the first entry of the last column makes the  columns of the $A-N$ block necessarily linearly independent, ensuring that this decomposition is in fact the sum of an invertible matrix and a nilpotent one.  The direct sum of such $m \times m$ blocks will respectively be invertible and nilpotent as well.

%%%%%%%%%%%%%%%%% Sharon's second part. Goes after general/non-zero blocks case %%%%%%%%%%%%%%%%%
Suppose the rational canonical form of the matrix is a series of companion matrix blocks followed by a zero block:
$$\begin{bmatrix}
           C_{g_1}  &   0     & \cdots        &  0\\
            0             &  \ddots &               & \vdots\\
            \vdots        &         & C_{g_m} & \vdots \\
            0             &   0     & 0             &[0]
      \end{bmatrix}$$

We will decompose the zero block of size $n \times n$ augmented with the last two rows and columns of the $C_{g_m}$ block of size $r \times r$ to ensure that there is at least one nonzero entry, which is the one on the last row of $C_{g_m}$ and the $(m-1)^{th}$ column.
Thus the augmented matrix that we then decompose is of the form
$$A = \begin{bmatrix}

	    0        &   -c_{s-2} & 0 & \cdots &  0\\
            1        &   -c_{s-1} &  0 & \cdots & 0 \\
            0        & \vdots       & \ddots & \ddots & \vdots \\
            \vdots & \vdots       & \ddots & \ddots & \vdots \\
	    0        &   0      & \cdots & \cdots & 0
      \end{bmatrix}$$
where A is an $(n+2) \times (n+2)$ matrix.

To find a suitable $(n+2) \times (n+2)$ nilpotent matrix of rank $n-1$ for the decomposition, we will conjugate $N_{n+2}$ by a suitable invertible matrix since conjugation will result in another nilpotent matrix.  Choose

 $$P = \begin{bmatrix}
            0      & \cdots & \cdots & \cdots &  0     & 1\\
            0      & 1      & -1     &  0      & \cdots & 0\\
            \vdots & 0      &  1     & -1     & \ddots & \vdots \\
            \vdots &  \vdots     & \ddots & \ddots & \ddots & 0\\
            0      &        &        & 0      & 1 & -1\\
            1      & -1     &  0      & \cdots & \cdots & 0
          \end{bmatrix} \text{ and } {P^{-1}}= \begin{bmatrix}
            1      & 1      & 1      & \cdots      & 1 & 1\\
            1      & 1      & 1      &        & \vdots & 0 \\
            1      & 0      & 1      & \ddots & \vdots & \vdots\\
            \vdots & \vdots &     \ddots  & \ddots & 1      & \vdots \\
            1      & \vdots      &        &  0      & 1      & \vdots\\
            1      & 0      & \cdots & \cdots & \cdots & 0
          \end{bmatrix}$$

Then,

 $$PN_{n+2}P^{-1} = \begin{bmatrix}
            1      & 0      & \cdots & \cdots & 0  & 1      & 0\\
            0      & 0      & \cdots & \cdots & 0  & 0      & 1\\
            0 & 1      &        &    & & 0 & 0 \\
            \vdots &        & \ddots &    & \vdots & \vdots\\
            \vdots &        &        &  \ddots&  & \vdots & \vdots \\
            0      &        &        &&    1& 0      & 0\\
            -1     & \cdots & \cdots & \cdots & -1 & -1     & -1
          \end{bmatrix}$$
 will be the nilpotent matrix involved in the decomposition.

Now let

$$U =
           \begin{bmatrix}
                       1             & -c_{s-2}     & 0 & \cdots &  0 &  1     & 0\\
                       1             & -c_{s-1}     & 0  & \cdots &  0 & 0      & 1\\
                       0             & 1    & 0 &   \cdots &  0      & 0 & 0 \\
                       \vdots        &        & \ddots &  &   & \vdots & \vdots\\
	               \vdots        &     &   & \ddots   &   & \vdots & \vdots\\
                       \vdots        &        & &     & & \vdots & \vdots \\
                       0             &        & &        &  1 & 0      & 0 \\
                       -1             & \cdots & \cdots   & \cdots & -1 & -1      & -1
                     \end{bmatrix} \text{ and } N=
                    \begin{bmatrix}
           -1      & 0      & \cdots & 0  & -1      & 0\\
            0      & 0      & \cdots & 0  & 0      & -1\\
            0 &  -1      &        &    & 0 & 0 \\
            \vdots &        & \ddots &    & \vdots & \vdots\\
            \vdots &        &        & -1  & \vdots & \vdots \\
            0      &        &        &    & 0      & 0\\
            1     & \cdots & \cdots & 1 & 1     & 1
          \end{bmatrix}$$

To see why U is invertible, note that regardless of the values of $c_{s-2}$ and $c_{s-1}$ the second through the last column form a linearly independent set of vectors because of the negative ones on different rows. Then we only need to determine if the first column can be written as a linear combination of vectors from this set. If indeed there were scalars such that the first column could be written as a linear combination of the others in $U$, then having zeros everywhere except the first, second and last row will restrict the coefficients of the columns to be zero except possibly the $(n+1)^{th}$  and $(n+2)^{th}$ column. However, a routine calculation shows such a linear combination is not possible either. Thus the matrix is invertible.

Now we consider A as the direct sum of $C_{g_1}$ through $C_{g_{(m-1)}}$ and $C_{g_m}$ with the zero block, as illustrated:
 $$\begin{bmatrix}
                        0  & \cdots      &  0 & -c_{0} & 0      \\
                       1  & \ddots &  \vdots &     \vdots & \vdots     \\
                           & \ddots & 0& -c_{s-2} & \vdots\\
                           &          & 1& -c_{s-1}  & \vdots\\

            \end{bmatrix} \oplus [0]$$
where $C_{g_m} \oplus [0]$ is an $(n+r) \times (n+r)$ matrix.

Then choose
 \setcounter{MaxMatrixCols}{20}
   $$U= C^*_{g_1} \oplus \cdots \oplus C^*_{g_{(m-1)}} \oplus  \begin{bmatrix}
               0  & \cdots     &  1 & 0 &  -c_{0}       \\
               1  & \ddots & \vdots & \vdots &     \vdots      \\
 	     		&\ddots   & 0 & \vdots & \vdots \\
                   &   & 1 & 1 & -c_{s-2} & 0  & \cdots &\cdots  & 0& 1& 0 \\
                           &  &        & 1& -c_{s-1} & 0&\cdots&\cdots &\cdots & 0&1\\
                    & & &    0 & 1    & 0 &   \cdots & \cdots   & 0      & 0 & 0 \\
                   & & &     \vdots        &        & \ddots &  &  & & \vdots & \vdots \\
	            & & &    \vdots        &     &   & \ddots   &  & & \vdots & \vdots \\
                    & & &    \vdots        &        & &        &\ddots  & & \vdots & \vdots \\
                    & & &    0             &        & &      &  &  1 & 0      & 0 \\
                      & & & -1             & \cdots & \cdots & \cdots  & \cdots & -1 & -1      & -1
            \end{bmatrix}$$

Then we have
$$N= N^* \oplus \cdots \oplus  \begin{bmatrix}
 0  & \cdots     &  -1 &        \\
                 & \ddots & \vdots      \\
 	     		&   & 0    \\
           & & & -1      & 0      & \cdots & 0  & -1      & 0\\
           & & &  0      & 0      & \cdots & 0  & 0      & -1\\
           & & &  0 &  -1      &        &    & 0 & 0 \\
           & & &  \vdots &        & \ddots &    & \vdots & \vdots\\
           & & &  \vdots &        &        & -1  & \vdots & \vdots \\
           & & &  0      &        &        &    & 0      & 0\\
           & & &  1     & \cdots & \cdots & 1 & 1     & 1
          \end{bmatrix}$$
which is nilpotent. Thus $A=U+N$ where $U$ is a unit and $N$ is a nilpotent.

One can check that in the last direct summand the first $(r-2)$ columns and the last $(n+1)$ columns are linearly independent because of the ones in different rows.  Due to the $1$'s the $r^{th}$ row and $(r+1)^{th}$ column, the $(r+1)^{th}$ column cannot be written as a linear combination of the others, thus the only concern is the $(r-1)$th column. Writing it as a linear combination of other columns requires constructing the linear combination using only the $(n+r-1)^{th}$ and $(r+n)^{th}$ columns, because to use the $(r-2)^{th}$ column requires also using the $r^{th}$ column due to the $1$ in the first row and $(r-2)^{th}$ column.  However, the $1$ in the $(r+1)^{th}$ row and $r^{th}$ column prevents this.  Yet we've observed the $(r-1)^{th}$, $(n+r-1)^{th}$, and $(r+n)^{th}$ columns are linearly independent. Thus $U$ is invertible.

%%%%%%%%%%%%%%%%%% Conclusion %%%%%%%%%%%%%%%%%

Having addressed all possible cases for an $n \times n$ matrix's rational canonical form, we again recall that conjugation by an invertible matrix and its inverse preserves both unity and nilpotence, so we may conclude that $\mathbb{M}_n(k)$ is nil-good for any field $k$ and dimension $n$.

\end{proof}
\end{theorem}
%%%%%%%%%%%%%%%%%% Part 3! %%%%%%%%%%%%%%%%%%%%

\section{Matrices over Division Rings}
We now consider linear operators on division ring modules of dimension $n$.  Since the modules we consider are finite-dimensional, there exists a basis with respect to which we may express linear operators as $n \times n$ matrices.

\begin{theorem}{For all $n \geq 1$ the ring $\mathbb{M}_n(D)$ is nil-good if $D$ is a division ring.}
\begin{proof}

Given an $n \times n$ matrix $A \in \mathbb{M}_n(D)$ there exists an invertible matrix $Q \in \mathbb{M}_n(D)$ such that $A=QA_dQ^{-1}$ where $A_d$ is a matrix of the form $ U_A \oplus N_A $.
 Here $U_A$ is an $m \times m$ invertible block on the diagonal and $N_A$ is an $(n-m) \times (n-m)$ nilpotent block.

Although a matrix over a division ring does not necessarily have a rational canonical form, there exists a primary rational canonical form \cite{Coh73} similar to the standard rational canonical form for certain matrices.  Suppose a matrix $A \in \mathbb{M}_n (R)$ is algebraic over the center of the division ring $R$ and that it has a single elementary divisor $\alpha$.  P. Cohn proved that if $\alpha = c_1c_2 \dots c_s$ then $A$ may be put into the form
$$\begin{bmatrix} C_1 & 0 &\cdots & 0 \\
 N^{*} & C_2 & \ddots & \vdots \\ 0 & \ddots & \ddots & 0\\ \vdots & \dots & N^{*} & C_s \end{bmatrix}$$
 where $N^{*}$ is a matrix with a $1$ in the entry in the upper right corner and zeros everywhere else, and $C$ is the companion matrix of $c_i$ for each $i$.

We note that the matrix $N_A$ may be put in primary rational ranonical form since it is algebraic over the center of any ring and its minimal polynomial $x^k = 0$ has a single elementary divisor.
%and we may express any $R$-module for a division ring $R$ as a direct sum of $R/({N_A}^k)$ %where the sum of the integers $k$ is the dimension $n$ since the quotient of $R$ by the ideal %generated by an invariant factor of $N_A$ is a submodule.
Since the companion matrix of any power of nilpotent $N_A$ is a direct sum of $N_r$ blocks, the primary rational canonical form of $N_A$ is zero everywhere except the subdiagonal, which, as in the case of matrices over a field, makes the primary rational canonical form a direct sum of $N_r$ blocks as well.

If the $N_r$ blocks have dimension $2$ or higher, then we can write $N_A$ in the form  $$\begin{bmatrix}
            N_{r_1}  & 0  &\cdots & 0       \\
            0 & N_{r_2}    & \ddots  & \vdots  \\
	    \vdots   &\ddots    & \ddots   &  0  \\
            0           & \cdots         &0    &N_{r_k}
          \end{bmatrix}$$
 which, by subtracting a nilpotent of the form
$$\begin{bmatrix}
            N^{*}_{r_1}  & 0  &\cdots & 0       \\
            0 & N^{*}_{r_2}    & \ddots  & \vdots  \\
	    \vdots   &\ddots    & \ddots   &  0  \\
            0           & \cdots         &0    &N^{*}_{r_k}
          \end{bmatrix}$$
 (in which the matrix $N^{*}_{r_i}$ is an $N^{*}$ matrix with the dimension of the corresponding $N_{r_i}$ matrix) results in an invertible matrix.  The details of this proccess are outlined more explicitly in the case of a vector space over a field.

If any companion matrix in the primary rational canonical form of $N_A$ has a $k \times k$ zero block in its companion matrix, then we consider the augmented diagonal block given by the zero block and the $k+1$ entries of the column immediately to the left and the $k+1$ entries furthest to the right in row immediately above.  If the zero block in the companion matrix is not the first diagonal block in the first companion matrix on the diagonal of the primary rational canonical form of $N_A$, then the augmented block in question will be of the form
$$\begin{bmatrix}
            0        & 0   & \cdots   & 0  \\
	     1  &0  & \cdots   &  0  \\
             0   &\vdots & \ddots    &\vdots    \\
            \vdots   &\ddots      &\cdots   &  0
          \end{bmatrix} $$
for which, if we substract the nilpotent matrix
$$ \begin{bmatrix}
           -1      & 0      & \cdots & 0  & -1      & 0\\
            0      & 0      & \cdots & 0  & 0      & -1\\
            0 &  -1      &        &    & 0 & 0 \\
            \vdots &        & \ddots &    & \vdots & \vdots\\
            \vdots &        &        & -1  & \vdots & \vdots \\
            0      &        &        &    & 0      & 0\\
            1     & \cdots & \cdots & 1 & 1     & 1
          \end{bmatrix}$$
the result is an invertible block.

Just as is the case for matrices over a field, this redistribution of columns and rows into slightly smaller or larger blocks does not change the block martices used to decompose $A_d$ into the sum of an invertible matrix and a nilpotent one.
In the case that the first $k \times k$ block of the first companion matrix on the diagonal of the primary rational canonical form of $N_A$ is a zero block, we treat the first $(m+k) \times (m+k)$ block, which is comprised of $U_A$ and the zero block in question, as described in the case detailed below, then treat the bottom $(n-m-k) \times (n-m-k)$ block as described above.

%%%%%%%%%%%%%%%%%%%%%%%%% Sharon's Part #3 %%%%%%%%%%%%%%%%%%%%%%%

In the special case that the $N_A$ block is an $(n-m) \times (n-m)$ zero block, a little more work is required than for the field case to find a decomposition of the matrix since the invertible block is not in a nice normal form.  We first consider the case when $a_{mm}$ is nonzero. Then we augment the zero block by the last $n-m$ entries of the last column and row of $U_{A}$.

The augmented matrix
$$A = \begin{bmatrix}
               a_{mm} &  0 & \cdots & 0 \\
            \vdots       & \ddots & \ddots & \vdots \\
            \vdots       & \ddots & \ddots & \vdots \\
	     0      & \cdots & \cdots & 0
      \end{bmatrix}$$
can be written as the sum of the following:\\

$$U =
           \begin{bmatrix}
                     a_{mm} +1     & 0  & \cdots & \cdots & 0 & 1      & 0\\
                     0        & 0 &   \cdots & \cdots   & 0      & 0 & 1 \\
                     \vdots & 1 &  &  & & \vdots & \vdots\\
	             \vdots &     & \ddots  &     & & \vdots & \vdots\\
                     \vdots &        & & \ddots      & & \vdots & \vdots \\
                      0       &        & &     &  1 & 0      & 0 \\
                       -1     & \cdots & \cdots  & \cdots & -1 & -1      & -1
                     \end{bmatrix} \text{ and } N=
                    \begin{bmatrix}
           -1      & 0      & \cdots & 0  & -1      & 0\\
            0      & 0      & \cdots & 0  & 0      & -1\\
            0 &  -1      &        &    & 0 & 0 \\
            \vdots &        & \ddots &    & \vdots & \vdots\\
            \vdots &        &        &   & \vdots & \vdots \\
            0      &        &        &   -1 & 0      & 0\\
            1     & \cdots & \cdots & 1 & 1     & 1
          \end{bmatrix}$$
In both matrices, the last $(n-m)$ columns of $U$ are linearly independent because of the ones on exclusively different rows. If writing the first column as a linear combination of the others is possible, only the $(n-m)$th column can be included. Since $a_{mm} \neq 0$, the first column is not a scalar multiple of the $(n-m)$th column, and thus $U$ is invertible.

Therefore,
$A_{d}=  U+N$, where

\setcounter{MaxMatrixCols}{20}
   $$U=  \begin{bmatrix}
               a_{11}  & \cdots     &  a_{1m}        \\
                & \ddots & \vdots       \\
 	     		&   & a_{mm}+1 & 0  & \cdots &\cdots  & 0& 1& 0\\
                      &        & 0&  0&\cdots&\cdots &\cdots & 0&1\\
                    & &    0 & 1    & 0 &   \cdots    & \cdots      & 0 & 0 \\
                   & &     \vdots                &  & \ddots & &  & \vdots & \vdots \\
	            & &    \vdots             &   &    & \ddots & & \vdots & \vdots \\
                     & &    0                     & &      &  &  1 & 0      & 0 \\
                      & & -1              & \cdots & \cdots  & \cdots & -1 & -1      & -1
            \end{bmatrix}$$

and
\setcounter{MaxMatrixCols}{20}
$$N=  \begin{bmatrix}
 0  & \cdots     &  0 &        \\
                 & \ddots & \vdots      \\
 	     		&   & 0    \\
           & & & -1      & 0      & \cdots & 0  & -1      & 0\\
           & & &  0      & 0      & \cdots & 0  & 0      & -1\\
           & & &  0 &  -1      &        &    & 0 & 0 \\
           & & &  \vdots &        & \ddots &    & \vdots & \vdots\\
           & & &  \vdots &        &        & -1  & \vdots & \vdots \\
           & & &  0      &        &        &    & 0      & 0\\
           & & &  1     & \cdots & \cdots & 1 & 1     & 1
          \end{bmatrix}.$$
 \\

One can see that the first $(m-1)$ and the last $n$ columns are linearly independent. So if $U$ is not invertible, then the $m^{th}$ column can be written as a linear combination of the other columns. Note that because of the zeros, the coefficients of the last $m$ columns, except possibly the $(m+1)^{th}$ column, have to be zero. Since there is a negative one on the last row of the $m^{th}$ column, the coefficient of the $(n-1)^{th}$ column has to be one. Then the difference of the $m^{th}$ column and the $(n-1)^{th}$ column would be a linear combination of the first $(m-1)$ columns. The existence of such a linear combination would imply that $U_{A}$ is not invertible, a contradiction. Therefore $U$ is invertible.

Now consider the case in which $a_{mm}=0$. If there is a nonzero entry on the diagonal, say $a_{ii}$, then we can conjugate $A_{d}$ by a permutation matrix $P$ that swaps the $i^{th}$ row with the $m^{th}$ row. Then the $(m,m)^{th}$ entry in $PA_{d}P^{-1}$ will be nonzero, and we can apply the above method to decompose the matrix.

If all of the entries on the diagonal are zero, but $a_{(m-1)m}\neq 0$, then we can conjugate $A_{d}$ by an invertible matrix $S$ that subtracts the $(m-1)^{th}$ row from the $m^{th}$ row. Note that conjugation does not change the invertibility of the unit block and the nilpotence of the nilpotent block.

We define $S$ and its inverse $S^{-1}$ as the following $(m+n) \times (m+n)$ matrices:
$$S = \begin{bmatrix}
1  \\
 & \ddots \\
&   & 1 \\
& & -1 & 1\\
& & & & I_n \\
\end{bmatrix} \text{ and } S^{-1} = \begin{bmatrix}
1  \\
 & \ddots \\
&   & 1 \\
& & 1 & 1\\
& & & & I_n \\
\end{bmatrix}$$

$$ \text{So } A_{d}' = SA_{d}S^{-1} =  \begin{bmatrix}
               0 & \cdots     &  a_{1(m-1)}+a_{1m}  &  a_{1m}      \\
                \vdots & \ddots & \vdots & \vdots      \\
 	     	a_{(m-1)1}	& \cdots  & a_{(m-1)m} &   a_{(m-1)m} \\
                a_{m1} - a_{(m-1)1}       & a_{m(m-1)}       & a_{m(m-1)}-a_{(m-1)m} &  -a_{(m-1)m}\\

                   \end{bmatrix} \oplus [0]$$

Now the $(m,m)^{th}$ entry of $U_{A_{d}'}$ is nonzero, so we can use the first case for decomposition.

If all of the entries on the diagonal are zero and but $a_{m-1(m)}= 0$, there is a nonzero entry on the $m$th column because $U_{A}$ is invertible. Let $a_{km}$ be nonzero. Then we can conjugate $A_{d}$ by a permutation matrix $P$ that swaps the $k^{th}$ row with the $(m-1)^{th}$ row.

$$ \text{Then let } A_{d}'' = PA_{d}P^{-1} =  \begin{bmatrix}
               0  & \cdots &a_{1(m-1)}   &\cdots   &  a_{1k}&a_{1m}      \\
                \vdots &  & \vdots & \vdots & \vdots & \vdots      \\
 	     	a_{(m-1)1}	& \cdots & 0& \cdots  & a_{(m-1)k} &   0 \\
                \vdots &  & \vdots & \vdots & \vdots & \vdots     \\
                a_{k1}	& \cdots  & a_{k(m-1)} & \cdots & 0 &  a_{km} \\
		a_{m1}     & \cdots       & a_{m(m-1)} & \cdots & a_{mk}&  0\\
                   \end{bmatrix} \oplus [0]$$

Now in the $m^{th}$ column, the entry on the $(m-1)^{th}$ row is nonzero, we can conjugate $A_{d}'$ by $S$ as introduced above and follow the above method for decomposition.

Therefore, solely by applying a certain series of invertible linear transformations, one may find a nil-good decomposition of any square matrix over a division ring.

\end{proof}
\end{theorem}

\section{General Properties of Nil-good Rings}

Having established the essential fact that the ring of $n \times n$ matrices over a division ring is nil-good, we may observe some sufficient or necessary conditions for some types of widely used rings to be nil-good.  In particular we give a necessary and sufficient condition for artinian rings and matrices over a local ring.  For completeness, proofs of other elementary facts are provided.

The following four remarks also appear in \cite{Dan15}, but we briefly provide our own proofs of these elementary results for completeness. 

\begin{remark}
A ring $R$ is nil-good if and only if there exists a nil ideal $\mathfrak{A}$ such that $R/ \mathfrak{A}$ is nil-good.
\begin{proof}
The forward direction is trivial, simply consider the ideal (0).
Suppose now that $\mathfrak{A}$ is a nil ideal of $R$ and every element of $R/ \mathfrak{A}$ has a nil-good decomposition.  If $\bar{a}$ is nilpotent in $R/\mathfrak{A}$ then $\bar{a}^k = 0$ in $R/\mathfrak{A}$ for some $k \in \mathbb{N}$ so $a^k \in R/\mathfrak{A}$.  If $R/\mathfrak{A}$ is a nil ideal, this implies $a^k=n$ for some nilpotent element $n$ in $R/\mathfrak{A}$.  Then if $a^k$ is nilpotent, it is immediate that $a$ is nilpotent.  Since any nil ideal is contained in the radical and units lift modulo $J(R)$, we conclude any unit $\bar{u}$ in $R/\mathfrak{A}$ lifts to a unit $u$ in $R$.  So the nil-good decomposition of any element in $R/\mathfrak{A}$ lifts to the sum of a unit and a nilpotent in $R$.

\end{proof}
\end{remark}

As a corollary to this, we know a ring $R$ is nil-good if and only if for any nil ideal $\mathfrak{A}$ the quotient ring $R/ \mathfrak{A}$ is nil-good.  The essence of the proof is similar to that of the above proposition, with the added observation that if $R$ is a nil-good ring and $\mathfrak{A}$ is a nil ideal, then given $a=n+u$ for a nilpotent $n$ in $R$ and a unit $u$ in $R$ the image $\bar{a}$ in $R/ \mathfrak{A}$ has the decomposition $\overline{u+n} = \bar{u} + \bar{n}$.  Since $\mathfrak{A}$ is a nil ideal it is contained in $J(R)$ so $\bar{u}$ is a unit in the quotient ring.  Moreover, the fact that $n^k = 0$ for some $k$ means that $n^k \in \mathfrak{A}$ so $\bar{n}^k = \bar{0}$.\\

\begin{remark}
If $R$ is nil-good then $J(R)$ is a nil ideal.
\begin{proof}
If $R$ is nil-good then for all $y \in J(R)$ we know $y=n+u$ where $n$ is nilpotent and $u$ is a unit or zero. Suppose for contradiction that $u$ is a unit.  Then if $U(R)$ denotes the set of units in $R$, we know $1-yu^{-1} \in U(R)$ by definition of the Jacobson radical. However $1-yu^{-1}  = 1 - nu^{-1} -1$ which implies $-nu^{-1} \in U(R)$, a contradiction. So we must have that $u=0$ which implies $y=n+0$ is nilpotent.
\end{proof}
\end{remark}

	If $J(R)$ is a nil ideal, then by the above two remarks we know $R$ is nil-good if and only if $R/J(R)$ is nil-good.  Therefore, if $J(R)$ is a nil ideal, we wish to know when $R/J(R)$ is nil-good. The following result will prove useful to that end.

If $J(R)$ is the unique ideal that is maximal both as a left ideal and as a right ideal, then we say $R$ is a local ring.\\

\begin{remark}
A local ring $R$ is nil-good if and only if $J(R)$ is nil.
\begin{proof}
If R is nil-good then $J(R)$ is nil by Remark 4.2.  If $J(R)$ is a nil ideal then since $J(R)$ is equal to the unique maximal ideal of $R$, any element not in $J(R)$ is a unit.  Therefore for any element $a$ in $R$ either has the decomposition $a=n+0$ or $a=u+0$ which implies $R$ is nil-good.
\end{proof}
\end{remark}

%If $J(R)$ is a nil ideal then we do know $R$ is nil-good if $J(R)$ is maximal both as a left ideal and as a right ideal.  If $J(R)$ is the unique ideal that is maximal both as a left ideal and as a right ideal, then we say $R$ is a local ring.  So as a corollary, we get that a local ring $R$ is nil-good if its unique maximal ideal is nil.  The converse is clearly true also, that if a local ring $R$ has a nil maximal ideal then $R$ is nil-good because all elements outside of the ideal are units.

\begin{remark}
If $R$ is nil-good then $R$ has no nontrivial central idempotents.
\begin{proof}
Suppose $R$ is nil-good and $e$ is a central idempotent. Then $e=u+n$ which implies $u=e-n$. If $u$ is a unit then $e-n$ is a unit that commutes with nilpotent $n$, which implies $u+n$ is a unit.  This implies $e=1$.  If $u=0$ then $e$ is nilpotent, but $e=e^2$ so $e=0$.  Therefore $e$ is trivial if it is central.
\end{proof}
\end{remark}

	The above remark allows us to conclude that a semisimple ring, which is always isomorphic to the direct product of matrix rings over division rings of various shapes and sizes, is nil-good if and only if it is simple.  Observing that if $R$ is a left artinian ring then $R/J(R)$ is semisimple \cite{Lam01}, we arrive at the following result.\\

 \begin{proposition}
If $R$ is a left artinian ring, then  $R$ is nil-good if and only if $J(R)$ is maximal.
\begin{proof}
Suppose $R$ is left artinian.  Then $R/J(R)$ is semisimple, and therefore $R \simeq \mathbb{M}_{n_1}(D_1) \times \dots \times \mathbb{M}_{n_r}(D_r)$ for some division rings $D_1,...,D_r$ and positive intgeres $n_1,...,n_r$.
However, any nontrivial direct product contains central idempotents.  So if $R$ is nil-good it must be that $R/J(R)$ is isomorphic to a matrix ring over a division ring, which is a simple semisimple ring.  The quotient $R/J(R)$ is simple if and only if the ideal $J(R)$ is maximal.
Conversely, $J(R)$ is nil if $R$ is artinian, and if $J(R)$ is maximal then we know $R/J(R)$ is a simple semisimple ring.  As shown for the first direction, we can conclude $R/J(R) \simeq \mathbb{M}_{n}(D)$ for some division ring $D$ and natural number $n$.  Then by Remark 4.1, the nil-good decompositions of $R/J(R)$ lift modulo $J(R)$. Thus $R$ is a nil-good ring.
\end{proof}
\end{proposition}

We may also establish a few facts about matrices over nil-good rings.  If $R$ is a simple artinian ring, then $\mathbb{M}_n(R)\simeq \mathbb{M}_n(\mathbb{M}_k(D))$ and since $\mathbb{M}_n(\mathbb{M}_k(D))\simeq \mathbb{M}_{nk}(D)$ we conclude that any matrix ring over a simple artinian ring is nil-good.

\begin{corollary}
If $R$ is a local ring such that $J(R)$ is nil, then $\mathbb{M}_n(R)$ is nil-good if and only if $\mathbb{M}_n(J(R))$ the maximal ideal of $\mathbb{M}_n(R)$ is a nil ideal.

\begin{proof}
Since $R$ is local $R/J(R)$ is a division ring, thus $\mathbb{M}_n(R/J(R))$ is nil-good.  Since ${\mathbb{M}_n(J(R))}$ the maximal ideal of $\mathbb{M}_n(R)$ is nil, by above remarks $\mathbb{M}_n(R)$ is nil-good.
Conversely, if $\mathbb{M}_n(R)$ is nil-good then by above proposition, its maximal ideal ${\mathbb{M}_n(J(R))}$ is nil.
\end{proof}

\end{corollary}

The above remark is included because although we would like to say that $R$ local and $J(R)$ nil implies that ${\mathbb{M}_n(R)}$ is nil-good, this only holds if the K\"{o}the conjecture does as well.

%%%%%%%%%%%%%%%%%%%%%%%%%%%%%%%%%%%%%%%%%%%%%%%%%%%%%%%%%%%%%%%%%%%

\section{The Nil-Good Clean Property}

%(Intro — Definition of NGC, explanation of Commutative Case)
We now define a new property ``nil-good clean,'' which is slightly weaker than clean or nil-good in general. An element is nil-good clean if it can be written as the sum of a nilpotent, idempotent and a unit. A ring is nil-good clean if all its elements are nil-good clean.

Observe that in the commutative case, nil-good clean is equivalent to the property clean. Let $R$ be a commutative ring that is nil-good clean. Then any element $a \in R$ can be written as the sum of a nilpotent $n$, an idempotent $e$ and a unit $u$. Since  $u' = u+n$ is always a unit in a commutative ring,
$a = u' + e$. Therefore, $a$ is clean. For the opposite direction, if $R$ is a clean commutative ring, then any element $a$ can be written as the sum of a unit $u$ and an idempotent $e$. It follows that by letting the nilpotent $n$ be zero, we have the nil-good clean decomposition $a = e+ u+ 0$.

%(Construction of LTDFM)
We have found one example of a nil-good clean ring that is not clean.  It is a subring of lower-triangular column-finite matrices.  We say a matrix is ``diagonal-finite'' if there is some fixed $n$ such that only the first $n$ sub-diagonals below the main diagonal contain nonzero entries for some nonnegative integer $n$.  The set of lower-triangular matrices with this property form a ring.  We denote the ring of column-finite matrices as $\mathbb{CFM_N}(R)$ and denote the element of this ring by $A=(a_{ij})^\infty_{i,j=1}$ for which $(a_{ij})_j$ are the rows, and $(a_{ij})_i$ are the columns.  \\

\begin{definition}
We denote the ring of lower-triangular diagonal-finite matrices over a ring $R$ by $\mathbb{LTDFM_N}(R) = \{A=(a_{ij})\in \mathbb{CFM_N}(R) |$ there exists $n\in \mathbb{N}$ such that $a_{ij}=0$ for all $i\geq j+n, j\geq 1\}$.
\end{definition}

To see that this is indeed a subring, note that if $A$ is a lower triangular matrix with $n$ nonzero subdiagonals and $B$ is a lower triangular matrix with $m$ nonzero subdiagonals then in the product $AB$, each column $(ab_{ij})_i$ will have $ab_{ij}=0$ above the diagonal, and $ab_{ij}=0$ if $i\geq j+n+m$ for all $j\geq 1$.  Therefore $AB$ can only have at most $n+m$ nonzero diagonals below the main diagonal, so the set is closed under multiplication. Clearly it is also a group with addition, and satisfies the usual ring axioms.\\

%(Propoition about idempotents and units in LTDFM)
\begin{lemma}For every idempotent $E$ in $\mathbb{LTDFM_{N}}(R)$ and $ i \in \mathbb{N}$, $e_{ii}$ is idempotent. Moreover, $ e_{ii} = 1$ for all $i$ implies $E = I$

\begin{proof}
As $E$ is idempotent, $E^2 = E$ and $e_{ii}e_{ii} = e_{ii}$. Thus $e_{ii}$ is idempotent.

We will show by induction on subdiagonls that $e_{ii}=1$ for all $i$ implies that the elements on the $n$th subdiagonal are zero. Consider the first subdiagonal as the base case. By the idempotence of $E$, when we consider the $(i-1)^{th}$ elment on the $i^{th}$ row, we have the equation $e_{i,i-1}(1)+(1)e_{i,i-1}=e_{i,i-1}$. Then $2e_{i,i-1} = e_{i,i-1}$  and $e_{i,i-1}=0$. This shows that all the elements on the first subdiagonal are zero. Suppose that all elements on or above the $k^{th}$ diagonal, except the main diagonal, are zero. We will now show that all elements on the $(k+1)^{th}$ diagonal are also zero, i.e. $e_{i,i-l}=0$ for any natural number $l \leq k$ and any $i$. Consider the $(i-k-1)^{th}$ element on the $i^{th}$ row. By the idempotence of $E$, $e_{i,i-k-1}=e_{i,i-k-1}e_{i-k-1,i-k-1}+e_{i,i-k}e_{i-k,i-k-1}+ \cdots +e_{i,i-1}e_{i-1,i-k-1}+e_{i,i}e_{i,i-k-1} = e_{i,i-k-1}(1)+(0)e_{i-k,i-k-1}+ \cdots + (0)e_{i-1,i-k-1}+(1)e_{i,i-k-1}=2 e_{i,i-k-1}$ and thus $e_{i,i-k-1}=0$. Therefore every element on the $(k+1)^{th}$ subdiagonal is zero. By induction, all the elements below the main diagonal are zero and $E = I$.
\end{proof}
\end{lemma}

\begin{lemma}For every unit $U$ in $\mathbb{LTDFM_{N}}(R)$ and $ i \in \mathbb{N}$, $u_{ii}$ is a unit.
\begin{proof}
Since U is invertible, there exists $A\in \mathbb{LTDFM_{N}}(R)$ such that $UA = I$. Then $u_{ii}a _{ii} = 1$ for any $ i \in \mathbb{N}$ and thus each $u_{ii}$ is a unit.
\end{proof}
\end{lemma}

\begin{theorem}
The ring $\mathbb{LTDFM_{N}}(R)$ is not exchange.
\begin{proof}
For contradiction suppose that $\mathbb{LTDFM_{N}}(R)$ is exchange. \\
Consider the matrix
$$ A = \begin{bmatrix} 1 &0 & \cdots & & & \\
 1 & 1 &0 &  \cdots  & &\\
0 & 1 & 1 & 0 & \cdots \\
\vdots & 0 & 1 & 1 & 0 & \cdots \\
 \vdots & \vdots & 0 & 1 & 1 & \ddots \\
\vdots & \vdots & \vdots & \ddots & \ddots & \ddots
\end{bmatrix}$$
We assume there exists an idempotent matrix $E$ such that $E \in ( \mathbb{LTDFM_{N}}(R) ) A$ and $I-E \in ( \mathbb{LTDFM_{N}}(R) ) (I-A)$.
Note that $I-E=B(I-A)$ for some $B \in  \mathbb{LTDFM_{N}}(R)$.

Since
$$I-A =  \begin{bmatrix} 0 &0 & \cdots & & & \\
 1 & 0 &0 &  \cdots  & &\\
0 & 1 & 0 & 0 & \cdots \\
\vdots & 0 & 1 & 0 & 0 & \cdots \\
 \vdots & \vdots & 0 & 1 & 0 & \ddots \\
\vdots & \vdots & \vdots & \ddots & \ddots & \ddots
\end{bmatrix}$$
 and $(I-A)_{ii}$ =0 for all $i$, $(I-E)_{ii} = B_{ii}((I-A)_{ii}) = B_{ii}(0) = 0$. Therefore $E_{ii} = 1$ for all $i$. As proved in the lemma earlier that $E$ with main diagonal of ones implies that $E=I$ and therefore  $I \in ( \mathbb{LTDFM_{N}}(R))A$. This means that there exists $C \in  \mathbb{LTDFM_{N}}(R)$ such that $CA = I$ and so A is invertible. The inverse of $A$ has to have ones and negative ones alternating infinitely in each column, and thus is not in $ \mathbb{LTDFM_{N}}(R)$. This is a contradiction so $ \mathbb{LTDFM_{N}}(R)$ is not exchange.
\end{proof}
\end{theorem}

%(RS is not Clean and an example)
\begin{corollary}{$\mathbb{LTDFM_{N}}(R)$ is not clean for any unital ring $R$}
\end{corollary}

%(Proof that LTDFM is NGC)
\begin{proposition}
The ring $\mathbb{LTDFM_N}(R)$ is nil-good clean if $R$ is a clean ring.
\begin{proof} Let $A=(a_{ij})\in \mathbb{LTDFM_N}(R)$ be an infinite matrix that has $n$ nonzero diagonals on or below the main diagonal.  We can write $A$ as the sum of two block-diagonal elements of $\mathbb{LTDFM_N}(R)$.  We define the first summand $D=\bigoplus D_k$.  It is the direct sum of $2n \times 2n$ blocks, thus for all $k\in \mathbb{N}$ define  $D_k=(d_{ij})^{2n(k+1)}_{i,j=1+2kn}$ where $d_{ij}=a_{ij}$ if $2kn<i,j\leq2(k+1)n$.  We define the second summand $N$ as the direct sum of the $n \times n$ zero matrix and $\bigoplus N_k$ where for all $k\in \mathbb{N}$ we define $N_k=(n_{ij})^{(2k+3)n}_{i,j=1+(2k+1)n}$ where $n_{ij}=a_{ij}$ if $2(k+1)n<i\leq(2k+3)n$ and $(2k+1)n<j<2(k+1)n$.

$$A =
\left[
\begin{array}{c@{}c@{}c}
  2n\left\{ \left[ \begin{array}{cccccccc}
         *& & & & & & &  \\
         *&*& & & & & &  \\
         *&*&*& & & & &  \\
         *&*&*&*& & & &  \\
         *&*&*&*&*& & &  \\
          &*&*&*&*&*& &  \\
          & &*&*&*&*&*&  \\
          & & &*&*&*&*&*  \\
  \end{array}\right] \right. &   &  \\
   2n \left\{ \left[ \begin{array}{cccccccccc}
          & & & & & &*&*&*&* \\
          & & & & & & &*&*&* \\
          & & & & & & & &*&* \\
          & & & & & & & & &* \\
          & & & & & & & & &  \\
           & & & & & & & & &  \\
           & & & & & & & & &  \\
           & & & & & & & & &  \\
           \end{array} \right] \right. &\left[ \begin{array}{cccccccc}
         *& & & & & & &  \\
         *&*& & & & & &  \\
         *&*&*& & & & &  \\
         *&*&*&*& & & &  \\
         *&*&*&*&*& & &  \\
          &*&*&*&*&*& &  \\
          & &*&*&*&*&*&  \\
          & & &*&*&*&*&*  \\
            \end{array}\right] &  \\
 &  & \ddots \\
\end{array} \right] $$

$$D =
\left[
\begin{array}{c@{}c@{}c}
  2n\left\{ \left[ \begin{array}{cccccccc}
         *& & & & & & &  \\
         *&*& & & & & &  \\
         *&*&*& & & & &  \\
         *&*&*&*& & & &  \\
         *&*&*&*&*& & &  \\
          &*&*&*&*&*& &  \\
          & &*&*&*&*&*&  \\
          & & &*&*&*&*&*  \\
  \end{array}\right] \right. &   &  \\
  &2n\left\{ \left[ \begin{array}{cccccccc}
         *& & & & & & &  \\
         *&*& & & & & &  \\
         *&*&*& & & & &  \\
         *&*&*&*& & & &  \\
         *&*&*&*&*& & &  \\
          &*&*&*&*&*& &  \\
          & &*&*&*&*&*&  \\
          & & &*&*&*&*&*  \\
            \end{array}\right] \right.  & \\
 &  & \ddots \\
\end{array}\right] $$

$$N =
\left[
\begin{array}{c@{}c@{}c@{}c}
  n\left\{ \left[ \begin{array}{ccc}
         & & \\
         &\mathbf{0}& \\
         & & \\
  \end{array}\right] \right. &   & & \\
   &2n\left\{ \left[ \begin{array}{cccccccc}
          & & & & & & &  \\
          & & & & & & &  \\
          & & & & & & &  \\
          & & & & & & &  \\
         *&*&*&*& & & &  \\
          &*&*&*& & & &  \\
          & &*&*& & & &  \\
          & & &*& & & &   \\ \end{array}\right] \right.  &  \\
  &  & 2n\left\{ \left[ \begin{array}{cccccccc}
          & & & & & & &  \\
          & & & & & & &  \\
          & & & & & & &  \\
          & & & & & & &  \\
         *&*&*&*& & & &  \\
          &*&*&*& & & &  \\
          & &*&*& & & &  \\
          & & &*& & & &   \\\end{array}\right] \right. &\\
          &&& \ddots\\
\end{array}\right] $$

Note that $N$ is a direct sum of finite nilpotent matrices, and therefore is nilpotent itself, and $D$ is the direct sum of finite matrices over a clean ring, so each $D_k$ has a clean decomposition, and therefore the direct sum of the units that decompose each $D_k$ is a unit, the direct sum of the idempotents that decompose each $D_k$ is an idempotent, and the sum of those two direct sums form a clean decomposition of $D$.  Therefore $\mathbb{LTDFM_N}(R)$ is nil-good clean.
\end{proof}
\end{proposition}

%Note that Bergman's Ring, and the example of a non-clean exchange ring that Ster derived from it, is not nil-good clean.  However, we can define a ring of endomorphisms that is analogous to $\mathbb{LTDFM_{N}}(R)$ over any ring whose elements can be represented as an infinite sequence of elements of a clean ring. The endomorphism ring will be nil-good clean for the same reason $\mathbb{LTDFM_{N}}(R)$ is.
%(Explanation)
%The nilpotent matrices, like the invertible ones and idempotent ones in that ring, must be upper triangular since the constant diagonals will not allow any nonzero subdiagonals.  So it is not nil-good clean for the same reason it is not clean.

\textbf{Question 1:} Is there a nil-good clean ring that is exchange but not clean? \\
We suspect that there is little overlap between nil-good clean rings that are not clean and exchange rings.\\

\textbf{Acknowledgement}

The authors would like to thank Alexander J. Diesl for his guidance and contributions to our research.\\

\nocite{*}
\bibliographystyle{amsplain}
\bibliography{math_references.bib}

%\textbf{References}\\

%\begin{tabular}{l p{0.8\linewidth}}
%[Coh73] & Paul M. Cohn, \textit{The similarity reduction of matrices over a skew field}, Math. Z. 132 (1973) 151-163.\\

%[Dan15] &  Peter Dancev, \textit{Nil-good unital rings}, pre-print, 2015.\\

%[Han77] & David Handelman, \textit{Perspectivity and ancellation in regular rings}, J. Algebra 48 (1977) 1-16.\\

%[Lam01] &  T.Y. Lam, \textit{A First Course in Noncommutative Rings}, Graduate Texts in Mathematics, Springer, 2001.\\

%[Nic77]  & W.K. Nicholson, \textit{Lifting idempotents and exchange rings}, Trans. Amer. Math. Soc. 229 (1977) 269-278.\\

%[Vam05]  & Peter Vamos, \textit{2-good rings}, Quart. J. Math. 56 (3) (2005) 417-430.\\

%[Wol53] & K. G. Wolfson, \textit{An ideal-theoretic characterization of the ring of all linear transformations}, Amer. J. Math. 75 (1953) 358–386. \\

%[Zel54] & D. Zelinsky, Every linear transformation is a sum of nonsingular ones, Proc. Amer. Math. Soc. 5 (1954) 627–630. \\

%\end{tabular}

\vspace{.3in}

%%%%%%%%%%%%%%%%%%%%%%%%%%%%%%%%%%%%%%%%%%%%%%%%%%%%%%%%%%%%%%%%%%%%

\end{document}